\documentclass[dvips]{amsart}
\usepackage{epsfig}

\newtheorem{thm}{Theorem}[section]
\newtheorem{cor}[thm]{Corollary}
\newtheorem{lem}[thm]{Lemma}

\theoremstyle{definition}

\theoremstyle{remark}
\newtheorem{rem}[thm]{Remark}
\numberwithin{equation}{section}

\newcommand{\N}{\mathbb{N}}

\newcommand{\C}{\mathbb{C}}
\newcommand{\ov}{\overline}
\newcommand{\dis}{\displaystyle}
\newcommand{\Notequiv}{/\kern-.6em\hbox{$\equiv$} }

\begin{document}

\title[Norms of products and factors polynomials]
{Norms of products and factors polynomials}%
\author{Igor E. Pritsker}%

\address{Department of Mathematics, 401 Mathematical Sciences, Oklahoma State
University, Stillwater, OK 74078-1058, U.S.A.}%
\email{igor@math.okstate.edu}

\thanks{Research supported in part by the National Science Foundation
grants DMS-9996410 and DMS-9707359.}%
\subjclass{Primary 11C08, 30C10; Secondary 30C85, 31A15}%
\keywords{Polynomials, products, factors, uniform norm, logarithmic capacity,
equilibrium measure, subharmonic function, Fekete points}%



\begin{abstract}

We study inequalities connecting a product of uniform norms of
polynomials with the norm of their product.  Generalizing
Gel'fond-Mahler inequality for the unit disk and
Kneser-Borwein inequality for the segment $[-1,1]$, we prove an
asymptotically sharp inequality for norms of products of
algebraic polynomials over an arbitrary compact set in plane.
Applying similar techniques, we produce a related inequality
for the norm of a single monic factor of a monic polynomial. The
best constants in both inequalities are obtained by potential
theoretic methods. We also consider applications of the general
results to the cases of a disk and a segment.

\end{abstract}

\maketitle


\section{Introduction}

Let $E$ be a compact set in the complex plane ${\C}$.  Define the uniform
(sup) norm on $E$ as follows:
$$\| f \|_E = \sup_{z \in E} |f(z)|.$$
Consider algebraic polynomials $\{ p_k (z) \}_{k =1}^m$ and their product
$$p(z) := \prod_{k =1}^m p_k (z).$$
We are interested here in polynomial inequalities of the form
\begin{equation} \label{1.1}
\prod_{k =1}^m \| p_k \|_E \leq C \| p \|_E .
\end{equation}
One of the first
results in this direction is due to Kneser \cite{Kn}, for $E = [-1,1]$ and
$m =2$ (see also Aumann \cite{Au}), who proved that
\begin{equation}  \label{1.2}
\| p_1 \|_{[-1,1]} \|p_2 \|_{[-1,1]} \leq K_{\ell ,n} \| p_1 p_2
\|_{[-1,1]},
\end{equation}
where
\begin{equation}  \label{1.3}
K_{\ell ,n} := 2^{n -1} \prod_{k =1}^{\ell} \left( 1 + \cos \frac{2k -1}{2n}
\pi \right) \prod_{k =1}^{n - \ell} \left( 1 + \cos \frac{2k -1}{2n} \pi
\right),
\end{equation}
$\deg p_1 = \ell$ and $\deg (p_1 p_2) =n$.  Note that (\ref{1.2}) becomes an equality
for the Chebyshev polynomial $t(z) = \cos n \arccos z = p_1 (z) p_2 (z)$, with
a proper choice of the factors $p_1 (z)$ and $p_2 (z)$. P. B. Borwein \cite{Bor}
gave a new proof of (\ref{1.2})-(\ref{1.3}) and generalized this to the
multifactor inequality
\begin{equation} \label{1.4}
\prod_{k =1}^m \| p_k \|_{[-1,1]} \leq 2^{n -1}
\prod_{k =1}^{[ \frac{n}{2} ]} \left( 1 + \cos \frac{2k -1}{2n} \pi \right)^2
\| p \|_{[-1,1]}.
\end{equation}
He has also showed that
\begin{equation} \label{1.5}
2^{n -1} \prod_{k =1}^{[ \frac{n}{2} ]} \left( 1 + \cos \frac{2k -1}{2n}
\pi \right)^2 \sim (3.20991 \ldots )^n, \mbox{ as } n \rightarrow \infty.
\end{equation}

Another case of the inequality (\ref{1.1}) was considered by Gel'fond
\cite[p. 135]{Ge} in
connection with the theory of transcendental numbers, for $E = \ov{D}$,
where $D := \{ w: |w| < 1 \}$ is the unit disk:
\begin{equation} \label{1.6}
\prod_{k =1}^m \| p_k \|_{\ov{D}} \leq e^n \| p \|_{\ov{D}} .
\end{equation}
The latter inequality was improved by Mahler \cite{Ma1},
who replaced $e$ by $2$:
\begin{equation} \label{1.7}
\prod_{k =1}^m \| p_k \|_{\ov{D}} \leq 2^n \| p \|_{\ov{D}} .
\end{equation}
It is easy to see that the base $2$ cannot be decreased, if $m =n$ and
$n \rightarrow \infty$.  However, (\ref{1.7}) has recently been further improved in
two directions. D. W. Boyd \cite{Boy1, Boy2} showed that, by taking in account the
number
of factors $m$ in (\ref{1.7}), one has
\begin{equation} \label{1.8}
\prod_{k =1}^m \| p_k \|_{\ov{D}} \leq (C_m)^n \| p \|_{\ov{D}},
\end{equation}
where
\begin{equation} \label{1.9}
C_m := \exp \left( \frac{m}{\pi} \int_0^{\pi/m} \log
\left(2 \cos \frac{t}{2}\right) dt \right)
\end{equation}
is asymptotically best possible for each fixed $m$, as $n \rightarrow \infty$.
Kro\'{o} and Pritsker \cite{KP}
showed that, for any $m \leq n,$
\begin{equation} \label{1.10}
\prod_{k =1}^m \| p_k \|_{\ov{D}} \leq 2^{n -1} \| p \|_{\ov{D}},
\end{equation}
where equality holds in (\ref{1.10}) for {\it each} $n \in {\N}$, with
$m =n$ and $p(z) = z^n -1$.
We give an asymptotically sharp inequality for the norm of products
of polynomials on arbitrary compact set in Section \ref{sec2}, which
generalizes the results of Mahler, Kneser and Borwein. This inequality and
other connected to it results were originally obtained in \cite{Pr1}.

A closely related problem is to estimate the norm of a single factor
via the norm of the whole polynomial. Clearly, we have to normalize the
problem by assuming that $p(z)$ is a monic polynomial of degree $n$, with
a monic factor $q(z)$, so that
$$ p(z)=q(z)\, r(z).$$
In the case of the unit disk, Boyd \cite{Boy1}
proved an asymptotically sharp inequality
\begin{equation} \label{1.11}
\| q \|_{\ov{D}} \leq \beta^n \| p \|_{\ov{D}},
\end{equation}
with
\begin{equation} \label{1.12}
\beta := \exp \left( \frac{1}{\pi} \int_0^{2\pi/3} \log \left(2
\cos \frac{t}{2}\right) dt \right).
\end{equation}
This inequality improved upon a series of results by Mignotte
\cite{Mi}, Granville \cite{Gr} and Glesser \cite{Gl}.

Further progress was made by Borwein in \cite{Bor}, for the
segment $[-a,a],\ a>0$ (see Theorems 2 and 5 there or see
Section 5.3 in \cite{BE}). In particular, Borwein proved that if
$\deg q = m$ then
\begin{equation} \label{1.13}
|q(-a)| \leq \| p \|_{[-a,a]} a^{m-n} 2^{n -1} \prod_{k =1}^{m}
\left(1+\cos \frac{2k -1}{2n} \pi \right),
\end{equation}
where the bound is attained for a monic Chebyshev polynomial of
degree $n$ on $[-a,a]$ and a factor $q$. He also showed that, for
$E=[-2,2]$, the constant in the above inequality satisfies
\begin{eqnarray*}
&&\limsup_{n \to \infty} \left(2^{m-1} \prod_{k =1}^{m} \left( 1
+ \cos \frac{2k -1}{2n} \pi \right)\right)^{1/n} \\ &\le& \lim_{n
\to \infty} \left(2^{[2n/3]-1} \prod_{k =1}^{[2n/3]} \left( 1 +
\cos \frac{2k -1}{2n} \pi \right)\right)^{1/n} \\ &=& \exp \left(
\int_0^{2/3} \log \left(2 + 2\cos{\pi x}\right) dx \right) =
1.9081 \ldots,
\end{eqnarray*}
which hints that
\begin{equation} \label{1.14}
C_{[-2,2]} = \exp \left( \int_0^{2/3} \log \left(2 + 2\cos{\pi
x}\right) dx \right) = 1.9081 \ldots.
\end{equation}
In Section \ref{sec3}, we find an asymptotically sharp inequality of this
type for a rather arbitrary compact set $E$. The general result is
then applied to the cases of a disk and a line segment, so that we
recover (\ref{1.11})-(\ref{1.12}) and confirm (\ref{1.14}). Also see
\cite{Pr2} for these results.

Considered problems have applications in transcendence theory (see \cite{Ge})
and in designing algorithms for factoring polynomials (see \cite{Boy3} and
\cite{La}). We confine ourselves to studying the sup norms for polynomials
of one variable only. A survey of the results involving other norms (e.g., Bombieri
norms) can be found in \cite{Boy3}. These inequalities are also of considerable
interest for polynomials in several variables, where very little is known
about sharp constants (cf. \cite{AM}, \cite{BBEM}, \cite{BST} and \cite{Ma2}).

\section{Products of Polynomials in Uniform Norms} \label{sec2}

Inequalities (\ref{1.2})-(\ref{1.10}) clearly indicate that the constant
$C$ in (\ref{1.1}) grows exponentially fast with $n$, with the base for
the exponential depending on the set $E$. A natural general problem arising
here is to find the {\it smallest} constant $M_E > 0,$ such that
\begin{equation} \label{2.1}
\prod_{k =1}^m \| p_k \|_E \leq M_E^n \| p \|_E
\end{equation}
for arbitrary algebraic polynomials $\{ p_k (z) \}_{k =1}^m$ with complex
coefficients, where $p(z) = \prod_{k =1}^m p_k (z)$ and $n = \deg p$.
The solution of this problem is based on the logarithmic potential theory
(cf. \cite{Ts} and \cite{Ra}). Let ${\rm cap}(E)$ be the {\it logarithmic
capacity} of a compact set $E \subset {\C}$. For $E$ with ${\rm cap}(E)>0$,
denote the {\it equilibrium measure} of $E$  by $\mu_E$. We remark that
$\mu_E$ is a positive unit Borel measure supported on $E$ (see \cite[p. 55]{Ts}).
Define
\begin{equation} \label{2.2}
d_E(z) := \max_{t \in E} |z -t|, \qquad z \in {\C},
\end{equation}
which is clearly a positive and continuous function on ${\C}$.

\begin{thm} \label{th2.1}
Let $E \subset {\C}$ be a compact set, ${\rm cap} (E) >0$.  Then the
best constant $M_E$ in {\rm (\ref{2.1})} is given by
\begin{equation} \label{2.3}
M_E = \frac{\exp\left(\int \log d_E(z) d \mu_E (z)\right)}{{\rm cap} (E)} .
\end{equation}
\end{thm}

One can see from (\ref{2.1}) or (\ref{2.3}) that $M_E$ is invariant with
respect to the rigid motions and dilations of the set $E$ in the plane.

Note that the restriction ${\rm cap} (E) >0$  excludes only very
{\it thin} sets from our consideration (see \cite[pp. 63-66]{Ts}), e.g.,
finite sets in the plane.  On the other hand, Theorem \ref{th2.1} is applicable to any
compact set with a connected component consisting of more than one point
(cf. \cite[p. 56]{Ts}). In particular, if $E$ is a continuum, i.e., a connected set,
then we obtain a simple universal bound for $M_E$.

\begin{cor} \label{cor2.1}
Let $E \subset \C$ be a bounded continuum (not a single point). Then we
have
\begin{equation} \label{2.4}
M_E \leq \frac{{\rm diam}(E)}{{\rm cap} (E)} \leq 4,
\end{equation}
where ${\rm diam}(E)$ is the Euclidean diameter of the set $E$.
\end{cor}

For the unit disk $D = \{ w: |w| < 1 \}$, we have that ${\rm cap}(\ov{D}) =1$
\cite[p. 84]{Ts} and that
\begin{equation} \label{2.5}
\mu_{\ov{D}} = \frac{1}{2 \pi} d \theta ,
\end{equation}
where $d \theta$ is the arclength on $\partial D$.  Thus
Theorem \ref{th2.1} yields
\begin{equation} \label{2.6}
M_{\ov{D}} = \exp\left(\frac{1}{2 \pi} \int_0^{2 \pi}
\log d_{\ov{D}}(e^{i \theta})\ d \theta\right) =
\exp\left(\frac{1}{2 \pi} \int_0^{2 \pi} \log 2\ d \theta\right) =2,
\end{equation}
so that we immediately obtain Mahler's inequality (\ref{1.7}).

If $E = [-1,1]$ then ${\rm cap}([-1,1]) = 1/2$ and
\begin{equation} \label{2.7}
\mu_{[-1,1]} = \frac{dx}{\pi \sqrt{1 -x^2}} , \quad x \in [-1,1],
\end{equation}
which is the Chebyshev (or arcsin) distribution (see \cite[p. 84]{Ts}).
Using Theorem \ref{th2.1}, we obtain
\begin{eqnarray} \label{2.8}
M_{[-1,1]} & = & 2\exp\left(\frac{1}{\pi} \int_{-1}^1
\frac{\log d_{[-1,1]}(x)}{\sqrt{1 -x^2}} dx\right) = 2\exp\left(\frac{2}{\pi}
\int_0^1 \frac{\log (1 +x)}{\sqrt{1 -x^2}} dx \right) \nonumber \\
& = & 2\exp\left(\frac{2}{\pi} \int_0^{\pi/2}
\log (1 + \sin t) dt \right) \approx 3.2099123 .
\end{eqnarray}
This gives the asymptotic version of Borwein's inequality (\ref{1.4})-(\ref{1.5}).

It appears that the upper bound 4 in Corollary \ref{cor2.1} is not the best possible.
One might conjecture that the sharp universal bounds are as follows
\begin{equation} \label{2.9}
2=M_{\ov D} \le M_E \le M_{[-1,1]} \approx 3.2099123,
\end{equation}
for any bounded non-degenerate continuum $E$.

It is of interest to determine the nature of the extremal polynomials for (\ref{2.1}).
We characterized the asymptotically extremal polynomials for (\ref{2.1}), i.e., those
polynomials, for which (\ref{2.1}) becomes an asymptotic equality as
$n \rightarrow \infty$, by their asymptotic zero distributions. The precise
statements of these results can be found in Theorems 2.3-2.5 and Corollaries
3.1-3.3 of \cite{Pr1}.

\section{Uniform Norm of a Single Factor} \label{sec3}

In the same way as in Section \ref{sec2}, we naturally arrive at the problem
to find the best (the smallest) constant $C_E$, such that
\begin{equation} \label{3.0}
\|q\|_{E} \le C_E^n \, \|p\|_E,\quad \deg p = n,
\end{equation}
is valid for {\it any} monic polynomial $p(z)$ and {\it any} monic factor $q(z)$.
Our solution of this problem is based on similar ideas, involving the
logarithmic capacity and the equilibrium measure of  $E$.

\begin{thm} \label{thm1}
Let $E \subset {\C}$ be a compact set, ${\rm cap} (E) >0$.  Then the best
constant $C_E$ in {\rm (\ref{3.0})} is given by
\begin{equation} \label{3.1}
C_E = \frac{\dis\max_{u \in \partial E} \exp\left( \int_{|z-u| \ge 1} \log
|z-u| d \mu_E (z)\right)}{{\rm cap} (E)}.
\end{equation}
Furthermore, if $E$ is regular then
\begin{equation} \label{3.2}
C_E = \max_{u \in \partial E} \exp\left(\dis -\int_{|z-u| \le 1} \log |z-u| d
\mu_E (z)\right).
\end{equation}
\end{thm}
The above notion of regularity is to be understood in the sense of exterior
Dirichlet problem (cf. \cite[p. 7]{Ts}).

One can readily see from (\ref{3.0}) or (\ref{3.1}) that the best constant
$C_E$ is invariant under the rigid motions of the set $E$ in the plane.
Therefore we consider applications of Theorem \ref{thm1} to the family of disks
$D_r := \{z:|z|<r\}$, which are centered at the origin, and to the family of
segments $[-a,a],\ a>0.$

\begin{figure}[htb]
  \centerline{\psfig{figure=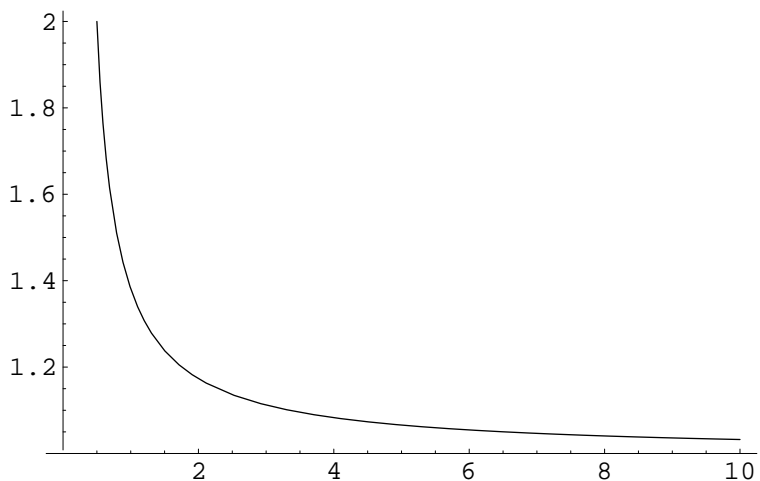,height=2.5in}}
  \caption{$C_{\ov D_r}$ as a function of $r$.}
  \label{fig1}
\end{figure}

\begin{cor} \label{cor1}
Let $D_r$ be a disk of radius $r$. Then the best constant $C_{\ov D_r}$, for
$E=\ov{D_r}$, is given by
\begin{equation} \label{3.3}
C_{\ov D_r} = \left \{ \begin{array}{l} \dis \frac{1}{r}, \quad 0<r \le 1/2, \\
\\ \dis \frac{1}{r} \exp \left(\frac{1}{\pi} \dis \int_0^{\pi-2\arcsin\frac{1}{2r}} \log
\left(2 r \cos\frac{x}{2}\right)\, dx \right), \quad  r>1/2.
\end{array} \right.
\end{equation}
\end{cor}
Note that (\ref{1.11})-(\ref{1.12}) immediately follow from (\ref{3.3}) for
$r=1.$ The graph of $C_{\ov D_r}$ is in Figure \ref{fig1}.

\begin{figure}[htb]
  \centerline{\psfig{figure=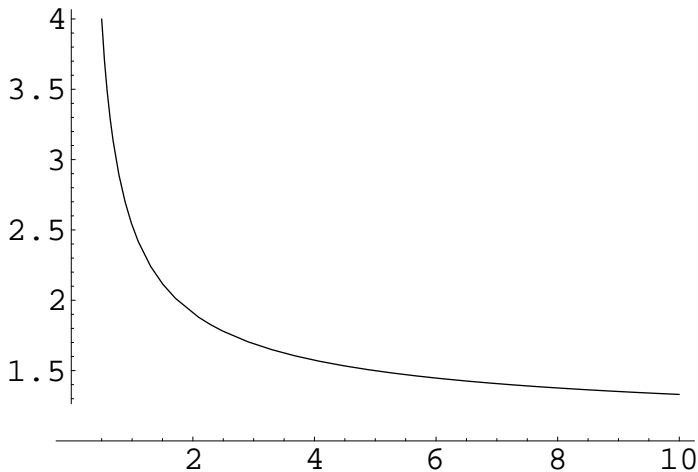,height=2.5in}}
  \caption{$C_{[-a,a]}$ as a function of $a$.}
  \label{fig2}
\end{figure}

\begin{cor} \label{cor2}
If $E=[-a,a],\ a>0$, then
\begin{equation} \label{3.4}
C_{[-a,a]} = \left \{ \begin{array}{l} \dis \frac{2}{a}, \quad 0< a \le 1/2,
\\ \\ \dis \frac{2}{a} \exp \left(\dis \int_{1-a}^a \frac{\log(t+a)}{\pi
\sqrt{a^2-t^2}}\, dt \right), \quad  a>1/2.
\end{array} \right.
\end{equation}
\end{cor}
Observe that (\ref{3.4}), with $a=2$, implies (\ref{1.14}) by the change of
variable $t=2\cos{\pi x}.$ We include the graph of $C_{[-a,a]}$ in Figure \ref{fig2}.

We now state two general consequences of Theorem \ref{thm1}. They explain
some interesting features of $C_E$, which the reader may have noticed in Corollaries
\ref{cor1} and \ref{cor2}. Recall that the Euclidean diameter of $E$ is defined by
$$ {\rm diam}(E) := \max_{z,\zeta \in E} |z-\zeta|.$$

\begin{cor} \label{cor3}
Suppose that ${\rm cap}(E)>0.$ If ${\rm diam}(E) \le 1$ then
\begin{equation} \label{3.5}
C_E = \frac{1}{{\rm cap}(E)}.
\end{equation}
\end{cor}
It is well known that cap$(D_r)=r$ and cap$([-a,a])=a/2$ (see \cite[p.
135]{Ra}), which clarifies the first lines of (\ref{3.3}) and (\ref{3.4}) by
(\ref{3.5}).

The next Corollary shows how the constant $C_E$ behaves under dilations of
the set $E$. Let $\alpha E$ be the dilation of $E$ with a factor $\alpha>0.$
\begin{cor} \label{cor4}
If $E$ is regular then
\begin{equation} \label{3.6}
\lim_{\alpha \to +\infty} C_{\alpha E} = 1.
\end{equation}
\end{cor}
Thus Figures \ref{fig1} and \ref{fig2} clearly illustrate (\ref{3.6}).

\medskip
We conclude this section with two remarks.

\begin{rem}
One can deduce inequalities of the type (\ref{3.0}), for various
$L_p$ norms, from Theorem \ref{thm1}, by using relations between $L_p$ and
$L_{\infty}$ norms of polynomials on $E$ (see, e.g., \cite{Pr3}).
\end{rem}

\begin{rem}
Note that the inequalities considered in this section hold for any monic factor
$q(z)$ of a monic polynomial $p(z)$, i.e., they hold for the {\it largest}
factor in the terminology of \cite{Boy4}. However, if we are granted the existence
of the factoring $p(z)=q(z)r(z)$, then the norm of the {\it smallest} factor (cf.
\cite{Boy5}) can be estimated from (\ref{2.1}) as follows:
\begin{equation} \label{3.7}
\|r\|_{E} \le M_E^{n/2} \, \|p\|_E^{1/2},\quad \deg p = n,
\end{equation}
which may be better than (\ref{3.0}) in some cases.
\end{rem}

\section{Proofs} \label{sec4}

The following lemma is a generalization of Lemma 2 in \cite{Boy2}.
We refer to \cite{Pr1} for its proof.

\begin{lem} \label{le5.1} Let $E \subset {\C}$ be a compact set
(not a single point) and let
$$d_E(z) := \max_{t \in E} |z-t|, \quad z \in \C.$$
Then $\log d_E(z)$ is a subharmonic function in ${\C}$ and
\begin{equation} \label{5.1}
\log d_E(z) = \int \log |z -t| d \sigma_E (t), \quad z \in {\C} ,
\end{equation}
where $\sigma_E$ is a positive unit Borel measure in ${\C}$ with unbounded
support, i.e.,
\begin{equation} \label{5.2}
\sigma_E({\C})=1 \quad \mbox{ and } \quad \infty \in {\rm supp}\, \sigma_E.
\end{equation}
\end{lem}

\begin{lem} \label{le5.2} {\rm (Bernstein-Walsh)}  Let $E \subset {\C}$
be a compact set, ${\rm cap}(E) >0$, with the unbounded component of
$\ov{{\C}} \setminus E$ denoted by $\Omega$.  Then, for any
polynomial $p(z)$ of degree $n$, we have
\begin{equation} \label{5.7}
|p (z)| \leq \| p \|_E \ e^{ng_{\Omega} (z, \infty )}, \quad
z \in {\C},
\end{equation}
where $g_{\Omega} (z, \infty )$ is the Green function of
$\Omega$, with pole at $\infty$, satisfying
\begin{equation} \label{5.8}
g_{\Omega} (z, \infty ) = \log \frac{1}{{\rm cap}(E)} + \int \log |z -t |
d \mu_E (t), \quad z \in {\C}.
\end{equation}
\end{lem}

This is a well known result about the upper bound (\ref{5.7}) for the growth
of $p(z)$ off the set $E$ (see \cite[p. 156]{Ra}, for example).  The
representation (\ref{5.8})
for $g_{\Omega} (z, \infty )$ is also classical (cf. Theorem III.37 in
\cite[p. 82]{Ts}).

Consider the $n$-th Fekete points $\{ a_{k,n} \}_{k =1}^n$ for a compact
set $E \subset {\C}$ (cf. \cite[p. 152]{Ra}).
Let
\begin{equation} \label{5.9}
F_n (z) := \prod_{k =1}^n (z - a_{k,n})
\end{equation}
be the Fekete polynomial of degree $n$, and define the normalized counting
measures in Fekete points by
\begin{equation} \label{5.10}
\tau_n := \frac{1}{n} \sum_{k =1}^n \delta_{a_{k,n}}, \quad n \in {\N},
\end{equation}
where $\delta_{a_{k,n}}$ is a unit point-mass at $a_{k,n}.$

\begin{lem} \label{le5.3} For a compact set $E \subset {\C}$,
${\rm cap} (E) >0$, we have that
\begin{equation} \label{5.11}
\lim_{n \rightarrow \infty} \| F_n \|_E^{1/n} = {\rm cap} (E)
\end{equation}
and
\begin{equation} \label{5.12}
\tau_n \stackrel{*}{\rightarrow} \mu_E, \quad \mbox{ as } n \rightarrow \infty.
\end{equation}
\end{lem}
Equation (\ref{5.11}) is standard (see Theorems 5.5.4 and 5.5.2 in \cite[pp. 153-155]{Ra}),
while (\ref{5.12}) follows from (\ref{5.11}) (see Ex. 5 on page 159 of \cite{Ra}).

\begin{proof}[Proof of Theorem \ref{th2.1}]
First we show that the best constant $M_E$ in (\ref{2.1}) is at most the
right hand side of (\ref{2.3}).  Clearly, it is sufficient to prove
an inequality of the type (\ref{2.1}) for monic polynomials only. Thus, we
assume that $p_k (z), 1 \leq k \leq m$, are all monic, so that $p(z)$ is monic too.
Let $\{z_{k,n}\}_{k=1}^n$ be the zeros of $p(z)$ and let
$\nu_n$ be the normalized zero counting measure for $p(z)$. Then,
we use (\ref{5.1}), Fubini's theorem and Lemma \ref{5.2} in the following estimate:
\begin{eqnarray*}
& &\frac{1}{n} \log \frac{\prod_{k=1}^m \| p_k \|_E}{\| p\|_E} \\ &\leq&
\frac{1}{n} \log \frac{\prod_{k=1}^n \| z - z_{k,n} \|_E}{\| p\|_E} =
\log \frac{1}{\| p\|_E^{1/n}} + \int \log d_E(z) d \nu_n(z)\\
&=& \log \frac{1}{\| p\|_E^{1/n}} + \int \int \log|z-t| d\nu_n(z) d\sigma_E(t)
= \int \log \frac{|p(t)|^{1/n}}{\| p\|_E^{1/n}} d \sigma_E(t)\\
&\leq& \int g_{\Omega}(t,\infty) d \sigma_E(t) =
\log \frac{1}{{\rm cap}(E)} + \int \int \log|z-t| d\sigma_E(t) d\mu_E(z)\\
&=& \log \frac{1}{{\rm cap}(E)} + \int \log d_E(z) d\mu_E(z).
\end{eqnarray*}
This gives that
\begin{equation} \label{5.18}
M_E \leq \frac{\exp ( \int \log d_E(z) d \mu_E (z))}{{\rm cap}(E)}.
\end{equation}

To show that equality holds in (\ref{5.18}), we consider the $n$-th Fekete points
$\{ a_{k,n} \}_{k =1}^n$ for $E$ and the Fekete polynomials $F_n(z), \ n \in {\N}.$
Observe that
$$\| z - a_{k,n} \|_E =d_E(a_{k,n} ), \quad 1 \leq k \leq n, \quad n \in {\N}.$$
Since ${\rm cap}(E) \neq 0$, the set $E$ consists of more than one point and, therefore,
$d_E(z)$ is a {\it strictly positive} continuous function in ${\C}$.
Consequently, $\log d_E(z)$ is also continuous in ${\C}$, and we obtain
by (\ref{5.12}) of Lemma \ref{le5.3} that
\begin{eqnarray} \label{5.18a}
& & \lim_{n \rightarrow \infty} \left( \prod_{k =1}^n \| z - a_{k,n} \|_E
\right)^{1/n} = \lim_{n \rightarrow \infty} \exp
\left( \frac{1}{n} \sum_{k =1}^n \log d_E(a_{k,n} ) \right)  \\ \nonumber
&=& \exp \left( \lim_{n \rightarrow \infty} \int \log d_E(z) d \tau_n (z) \right)
= \exp \left( \int \log d_E(z) d \mu_E (z) \right) .
\end{eqnarray}
Finally, we have from the above and (\ref{5.11}) that
$$M_E \geq \lim_{n \rightarrow \infty} \frac{( \prod_{k =1}^n \|
z -a_{k,n} \|_E )^{1/n}}{\| F_n \|_E^{1/n}} =
\frac{\exp ( \int \log d_E(z)d \mu_E (z))}{{\rm cap}(E)}.$$
\end{proof}

\begin{proof}[Proof of Corollary \ref{cor2.1}]
Since $E$ is a bounded continuum, we obtain from Theorem 5.3.2(a) of
\cite[p. 138]{Ra} that
$$ {\rm cap}(E) \geq \frac{{\rm diam}(E)}{4}.$$
Thus, the Corollary follows by combining this estimate with the obvious
inequality
$$ d_E(z) \leq {\rm diam}(E), \quad z \in E,$$
and by using that $\mu_E({\C})=1,\ {\rm supp}\mu_E \subset E.$
\end{proof}

\begin{proof}[Proof of Theorem \ref{thm1}]
The proof of this result is quite similar to that of Theorem \ref{th2.1} (also see \cite{Boy1}).
For $u \in \C$, consider a function
$$ \rho_u(z):=\max (|z-u|,1), \quad z \in \C.$$
One can immediately see that $\log \rho_u(z)$ is a subharmonic function in $z \in \C$,
which has the following integral representation (see \cite[p. 29]{Ra}):
\begin{equation} \label{4.1}
\log \rho_u(z) = \int \log |z-t|\, d\lambda_u(t), \quad z \in \C,
\end{equation}
where $d\lambda_u(u+e^{i\theta})=d\theta/(2\pi)$ is the normalized angular
measure on $|t-u|=1$.

Let $u \in \partial E$ be such that
$$ \|q\|_E=|q(u)|.$$
If $z_k,\ k =1, \ldots ,n,$ are the zeros of $p(z)$, arranged so that the first $m$ zeros
belong to $q(z)$, then
\begin{eqnarray} \label{4.2}
\log \|q\|_E &=& \sum_{k =1}^m  \log |u-z_k| \leq \sum_{k =1}^m \log \rho_u(z_k)
\leq \sum_{k =1}^n \log \rho_u(z_k)\nonumber \\
&=& \sum_{k =1}^n \int \log |z_k -t|\, d\lambda_u(t) = \int \log |p(t)|\, d\lambda_u(t),
\end{eqnarray}
by (\ref{4.1}).

It follows from (\ref{4.1})-(\ref{4.2}), Lemma \ref{le5.2} and Fubini's theorem that
\begin{eqnarray*}
\frac{1}{n} \log \frac{\|q\|_E}{\| p\|_E} &\leq&
\int \log \frac{|p(t)|^{1/n}}{\| p\|_E^{1/n}}\, d \lambda_u(t)
\leq \int g_{\Omega}(t,\infty)\, d \lambda_u(t) \\ &=&
\log \frac{1}{{\rm cap}(E)} + \int \int \log|z-t|\, d\lambda_u(t) d\mu_E(z)\\
&=& \log \frac{1}{{\rm cap}(E)} + \int \log \rho_u(z)\, d\mu_E(z).
\end{eqnarray*}
Using the definition of $\rho_u(z)$, we obtain from the above estimate that
\begin{eqnarray*}
\| q \|_E &\leq& \left( \frac{\dis \max_{u \in \partial E}\exp\left(\int\log\rho_u(z)\, d\mu_E (z)\right)}
{{\rm cap}(E)} \right)^n \| p \|_E \\ &=& \left(\frac{\dis\max_{u \in \partial E} \exp\left(
 \int_{|z-u| \ge 1} \log |z-u|\, d \mu_E (z)\right)}{{\rm cap} (E)}\right)^n \| p \|_E.
\end{eqnarray*}
Hence
\begin{equation} \label{4.5}
C_E \le \frac{\dis\max_{u \in \partial E} \exp\left( \int_{|z-u| \ge 1} \log
|z-u|\, d \mu_E (z)\right)}{{\rm cap} (E)}.
\end{equation}

In order to prove the inequality opposite to (\ref{4.5}), we consider the
$n$-th Fekete points $\{ a_{k,n} \}_{k =1}^n$ for the set $E$ and the Fekete
polynomials $F_n(z),\ n \in \N$.
Let $u \in \partial E$ be a point, where the maximum of the right hand side of
(\ref{4.5}) is attained. Define the factor $q_n(z)$ for $F_n(z)$, with zeros being
the $n$-th Fekete points satisfying $|a_{k,n}-u| \ge 1$. Then we have by
(\ref{5.12}) that
\begin{eqnarray*}
\lim_{n \rightarrow \infty} \| q_n \|_E^{1/n} &\ge& \lim_{n \rightarrow \infty} |q_n(u)|^{1/n} =
\lim_{n \rightarrow \infty} \exp\left( \frac{1}{n} \sum_{|a_{k,n}-u| \ge 1}\log|u-a_{k,n}| \right)  \\
&=& \exp \left( \lim_{n \rightarrow \infty} \int_{|z-u| \ge 1} \log|u-z|\, d \tau_n (z)
\right) \\ &=& \exp \left( \int_{|z-u| \ge 1} \log|u-z| d \mu_E (z) \right).
\end{eqnarray*}
Combining the above inequality with (\ref{5.11}) and the definition of $C_E$, we
obtain that $$C_E \geq \lim_{n \rightarrow \infty} \frac{\|q_n\|_E^{1/n}}{\|
F_n \|_E^{1/n}} \ge \frac{\dis \exp\left(\int_{|z-u| \ge 1} \log |z-u|\, d
\mu_E (z)\right)}{{\rm cap} (E)}.$$ This shows that (\ref{3.1}) holds true.
Moreover, if $u \in \partial E$ is a regular point for $\Omega$, then we obtain
by Theorem III.36 of \cite[p. 82]{Ts}) and (\ref{5.8}) that $$ \log
\frac{1}{{\rm cap}(E)}+\int \log|u-t|\, d\mu_E(t) = g_{\Omega}(u,\infty) = 0.$$
Hence $$ \log \frac{1}{{\rm cap}(E)}+\int_{|z-u| \ge 1} \log|u-t|\, d\mu_E(t) =
-\int_{|z-u| \le 1} \log|u-t|\, d\mu_E(t), $$ which implies (\ref{3.2}) by
(\ref{3.1}).
\end{proof}

\begin{proof}[Proof of Corollary \ref{cor1}]
It is well known \cite[p. 84]{Ts} that cap$(\ov {D_r})=r$ and $d \mu_{\ov
{D_r}}(re^{i\theta}) = d\theta/(2\pi),$ where $d\theta$ is the angular measure on $\partial
{D_r}$. If $r \in (0,1/2]$ then the numerator of (\ref{3.1}) is equal to 1, so
that
$$ C_{\ov {D_r}} = \frac{1}{r}, \quad 0<r\le1/2. $$
Assume that $r>1/2.$ We set $z=re^{i\theta}$ and let $u_0=re^{i\theta_0}$ be a
point where the maximum in (\ref{3.1}) is attained. On writing
$$ |z-u_0|=2r\left|\sin\frac{\theta-\theta_0}{2}\right|,$$
we obtain that
\begin{eqnarray*}
C_{\ov {D_r}} &=& \frac{1}{r}\exp\left(\frac{1}{2\pi}\int_{\theta_0+2\arcsin\frac{1}{2r}}
^{2\pi+\theta_0-2\arcsin\frac{1}{2r}}
\log\left|2r\sin\frac{\theta-\theta_0}{2}\right|\, d\theta\right) \\ &=&
\frac{1}{r} \exp \left(\frac{1}{2\pi} \dis \int_{2\arcsin\frac{1}{2r}-\pi}^{\pi-2\arcsin\frac{1}{2r}}
\log \left(2 r \cos\frac{x}{2}\right) dx \right) \\ &=&
\frac{1}{r} \exp \left(\frac{1}{\pi} \dis \int_0^{\pi-2\arcsin\frac{1}{2r}} \log
\left(2 r \cos\frac{x}{2}\right) dx \right),
\end{eqnarray*}
by the change of variable $\theta-\theta_0=\pi-x.$
\end{proof}

\begin{proof}[Proof of Corollary \ref{cor2}]
Recall that cap$([-a,a])=a/2$ (see \cite[p. 84]{Ts}) and
$$ d\, \mu_{[-a,a]}(t)=\frac{dt}{\pi\sqrt{a^2-t^2}}, \quad t \in [-a,a].$$
It follows from (\ref{3.1}) that
\begin{equation} \label{4.8}
C_{[-a,a]} = \frac{2}{a}\, \exp \left(\dis\max_{u \in [-a,a]} \int_{[-a,a]\setminus(u-1,u+1)}
\frac{\log|t-u|} {\pi \sqrt{a^2-t^2}}\, dt \right).
\end{equation}
If $a \in (0,1/2]$ then the integral in (\ref{4.8}) obviously vanishes, so that
$C_{[-a,a]} = 2/a$. For $a>1/2$, let
\begin{equation} \label{4.9}
f(u):= \int_{[-a,a]\setminus(u-1,u+1)} \frac{\log|t-u|} {\pi \sqrt{a^2-t^2}}\, dt.
\end{equation}
One can easily see from (\ref{4.9}) that
$$ f'(u) = \int_{u+1}^a \frac{dt} {\pi(u-t) \sqrt{a^2-t^2}} < 0, \quad u \in
[-a,1-a],$$
and
$$ f'(u) = \int_{-a}^{u-1} \frac{dt} {\pi(u-t) \sqrt{a^2-t^2}} > 0, \quad u \in
[a-1,a].$$
However, if $u \in (1-a,a-1)$ then
$$f'(u) = \int_{u+1}^a \frac{dt} {\pi(u-t) \sqrt{a^2-t^2}} +
\int_{-a}^{u-1} \frac{dt} {\pi(u-t) \sqrt{a^2-t^2}}.$$
It is not difficult to verify directly that
$$\int \frac{dt} {\pi(u-t) \sqrt{a^2-t^2}} = \frac{1}{\pi\sqrt{a^2-u^2}} \log
\left|\frac{a^2-ut+\sqrt{a^2-t^2}\,\sqrt{a^2-u^2}}{t-u}\right| + C, $$
which implies that
$$ f'(u)= \frac{1}{\pi\sqrt{a^2-u^2}} \log \left( \frac{a^2-u^2+u+\sqrt{a^2-(u-1)^2}\,\sqrt{a^2-u^2}}
{a^2-u^2-u+\sqrt{a^2-(u+1)^2}\,\sqrt{a^2-u^2}} \right),  $$
for $u \in (1-a,a-1).$ Hence
$$ f'(u)<0,\ u \in (1-a,0), \quad \mbox{ and } \quad f'(u)>0,\ u \in (0,a-1). $$
Collecting all facts, we obtain that the maximum for $f(u)$ on $[-a,a]$ is
attained at the endpoints $u=a$ and $u=-a$, and it is equal to
$$\max_{u \in [-a,a]} f(u) = \int_{1-a}^a \frac{\log(t+a)}{\pi \sqrt{a^2-t^2}}\, dt.$$
Thus (\ref{3.4}) follows from (\ref{4.8}) and the above equation.
\end{proof}

\begin{proof}[Proof of Corollary \ref{cor3}]
Note that the numerator of (\ref{3.1}) is equal to 1, because $|z-u| \le 1,\
\forall\, z \in E,\ \forall\, u \in \partial E$. Thus (\ref{3.5}) follows
immediately.
\end{proof}

\begin{proof}[Proof of Corollary \ref{cor4}]
Observe that $C_E \ge 1$ for any $E \in \C$, so that $C_{\alpha E} \ge 1$.
Since $E$ is regular, we use the representation for $C_E$ in (\ref{3.2}). Let
$T:E \to \alpha E$ be the dilation mapping. Then $|Tz-Tu|=\alpha|z-u|,\ z,u \in E,$
and $d\mu_{\alpha E}(Tz)=d\mu_E(z)$. This gives that
\begin{eqnarray*}
C_{\alpha E} &=& \max_{Tu \in \partial (\alpha E)} \exp\left(\dis -\int_{|Tz-Tu| \le 1}
\log |Tz-Tu|\, d\mu_{\alpha E} (Tz)\right) \\ &=&
\max_{u \in \partial E} \exp\left(\dis -\int_{|z-u| \le 1/\alpha}
\log (\alpha |z-u|)\, d\mu_{E} (z)\right) \\ &=&
\max_{u \in \partial E} \exp\left(\dis -\mu_E(\ov{D_{1/\alpha}(u)})\log\alpha-
\int_{|z-u| \le 1/\alpha} \log |z-u|\, d\mu_{E} (z)\right) \\ &<&
\max_{u \in \partial E} \exp\left(\dis -\int_{|z-u| \le 1/\alpha} \log |z-u|\, d\mu_{E}
(z)\right),
\end{eqnarray*}
where $\alpha \ge 1$. Using the
absolute continuity of the integral, we have that
$$ \lim_{\alpha \to +\infty} \int_{|z-u| \le 1/\alpha} \log |z-u|\, d\mu_{E} (z)
=0, $$
which implies (\ref{3.6}).
\end{proof}


\end{document}